\documentclass[12pt,noinfoline]{imsart}

\usepackage{amscd}
\usepackage{amsmath, amssymb,amsthm, amscd, amsfonts}
\usepackage{graphicx, caption, subcaption}
\usepackage{geometry}
\usepackage{url, hyperref, color}
\usepackage{here}
\usepackage{verbatim}

\theoremstyle{plain} 
\newtheorem{thm}{Theorem}[section]

\theoremstyle{definition}

\newtheorem{alg}[thm]{Algorithm}

\newtheorem{example}[thm]{Example}

\newcommand{\nn}[0]{\hspace*{.7em}}

\begin{document}

\pagestyle{plain}
\title{$\beta$ models for random hypergraphs with a given degree sequence}
\author{Despina Stasi\thanks{\textit{despina.stasi@gmail.com}, Illinois Institute of Technology},
Kayvan Sadeghi\thanks{\textit{kayvans@andrew.cmu.edu}, Carnegie Mellon University},
Alessandro Rinaldo\thanks{\textit{arinaldo@stat.cmu.edu}, Carnegie Mellon University},
Sonja Petrovi\'{c}\thanks{\textit{sonja.petrovic@iit.edu}, Illinois Institute of Technology},
Stephen E. Fienberg\thanks{\textit{fienberg@stat.cmu.edu}, Carnegie Mellon University}
}


\begin{abstract}
We introduce the 
beta model for random hypergraphs in order to
represent
the occurrence of multi-way interactions among agents in a social network. This
model builds upon and generalizes the well-studied beta model for random graphs,
which instead only considers pairwise interactions. 
We provide two algorithms for fitting the model parameters, IPS (iterative proportional scaling) and fixed point algorithm, 
prove that both algorithms converge if maximum likelihood estimator (MLE) exists, and provide algorithmic and geometric ways of dealing the issue of MLE existence.
\end{abstract}


\maketitle

\section{Introduction}
\begin{figure}[h]
\centering
\includegraphics[scale=1]{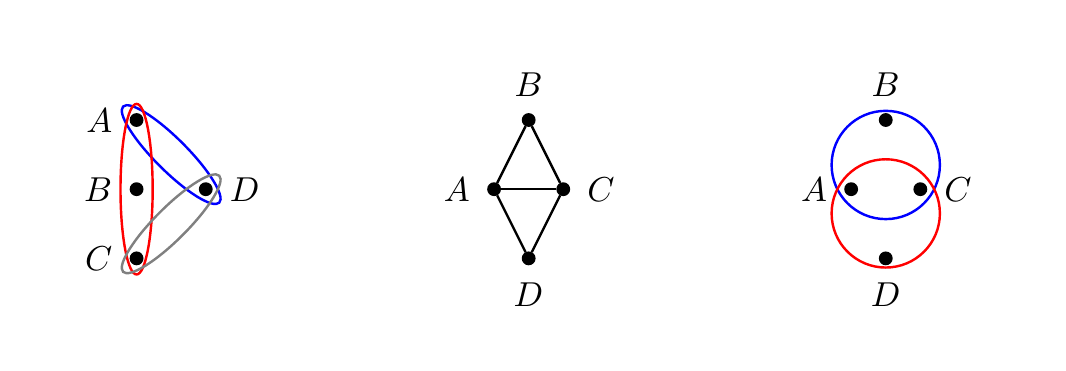}
\caption{Distinct hypergraphs $H$ and $H'$ reduced to  same graph $G$ (left, right, middle). 
}
\label{fig1}
\end{figure}

Social network models \cite{kolBook} are statistical
models for the joint occurrence of random edges in a graph, as a means to model
social interactions among agents in a population of interest. These
models typically focus on representing only \textit{binary} relations between
individuals.
As a result, when one is interested in higher-order ($k$-ary) interactions,
statistical models based on graphs may be ineffective or inadequate. 
Examples of $k$-ary relations are plentiful, and include forum or committee membership, co-authorship on scientific papers, or proximity of groups of people in photographs.
These datasets have  been studied by replacing each $k$-dimensional group with a number of binary relations (in particular, $k\choose 2$ of them, which form a clique), thus extracting binary information from the data, and then modeling and studying the resulting graph. Such a process  inevitably causes information loss.
For instance, let us consider statisticians Adam ($A$), Barbara ($B$), Cassandra
($C$), and David ($D$), see Figure~\ref{fig1}. Suppose the authors wrote three papers in following groups: 
  $(A,B,C)$, $(A,D)$, $(C,D)$. Representing this
information as a graph with edges between any two individuals who have co-authored
a paper  provides a graph $G$ with edges $\{(A,B),
(B,C), (C,D), (A,C)\}$. A hypergraph $H$ representing this information would
instead use the exact groups as hyperedges and, unlike $G$, would
be able to represent additional properties of such interactions, including how many papers were coauthored by these four individuals;  see Figure~\ref{fig1}.
If, in addition, $A$ is more likely to write a 3-author paper than a 2-author paper, 
this requires modeling separately the probabilities of these collaborations.
Despite the growing needs of practical values, models for random hypergraphs are relatively few
and simple. Random hypergraphs have been studied (\cite{KaronskiLuczak02}) as generalizations of the simple Erd\"os-R\'enyi model 
\cite{erd59} for networks; \cite{ghoshal2009random} considers an application of random tripartite hypergraphs to Flickr photo-tag data.

In this paper we introduce a simple and natural class of statistical models for
random hypergraphs, which we term hypergraph beta models, that allows one to model directly simultaneous higher-order
(and not only binary) interactions  among individuals in a network. As its name suggests, our model arises as a natural extension of the well-studied beta model for random graphs,  the exponential family
 for undirected networks which assumes independent edges and
whose minimal sufficient statistics vector is the degree sequence of the graph.
It is a special class of the more general  of $p_1$ models 
 \cite{hol81} which
assume independent edges and parametrize the probability of each edge by the
propensity of the two endpoint nodes. This model has been studied extensively; see \cite{bli09,ChatterjeeDS2011BetaModel,rin13,YanXuYang12,YanXu12}, which give, among other results, methods for model fitting. 

 Below we formalize the class of the beta models for hypergraphs. 
Just like the graph beta model, these are natural exponential
random graph models over hypergraphs which postulate independent edges and
whose sufficient statistics are the (hypergraph) degree sequences.
Our contributions are two-fold: first we formalize three classes of linear
exponential families  for
random hypergraphs of increasing degree of complexity and derive the
corresponding sufficient statistics and moment equations for obtaining the maximum likelihood estimator (MLE)
of the model parameters. Secondly, we design two iterative algorithms for
fitting these models that do not require evaluating the gradient or Hessian of the likelihood
function and can therefore be applied to large data: a variant of the IPS algorithm and a fixed point iterative algorithm  to compute the MLE of the edge probabilities and of the natural parameters, respectively. 
 We show that both algorithms will converge if the MLE exists. Finally, we
 illustrate our results and methods with some simulations.

As our analysis reveals, the study of the theoretical and asymptotic properties of hypergraph beta models
is especially challenging, more so than with the ordinary beta model.
The complexity of the new models, in turn,  leads to the problem of optimizing a
 complex likelihood function. Indeed, when the MLE does not exist,
optimizing the
likelihood function becomes highly non-trivial and, to a large extent, unsolved
for our model as well as for many other discrete linear exponential families.
To this end, we describe a geometric way  for dealing with the issue of existence of the MLE for these models and gain further insights into this difficult problem with simulation experiments.

\section{The hypergraph beta model: three variants}\label{sec:beta-models}

A \emph{hypergraph} $H$ is a pair $(V,F)$, where $V=\{v_1,\ldots,v_n\}$
is a set of \emph{nodes} (vertices) and $F$ is a family of non-empty subsets
of $V$ of cardinality different than $1$; the elements of $F$ are called the \emph{hyperedges} (or simply \emph{edges})
of $H$.   In a \emph{$k$-uniform} hypergraph, all edges are of size $k$. 
We restrict ourselves to the set $\mathcal{H}_n$  of hypergraphs on $n$ nodes, where nodes have a distinctive labeling. 
  Let $E = E_n$ be the set of all realizable hyperedges for a
hypergraph on $n$ nodes. While $E$ can in principle be the set of all possible
hyperedges, below we will consider more parsimonious models in which $E$ is
restricted to be a structured subset of edges. Thus we may write a hypergraph $x= (V,F)\in \mathcal{H}_n$ as the zero/one vector $x = \{x_e, e \in E \}$, where $x_e=1$
for $e\in F$ and $x_e=0$ for $e\in E\setminus F$.
The degree of a node in $x$ is the number of edges it belongs to; the degree information for $x$ is summarized in the degree sequence vector whose $i$th entry is the degree $d_i(x)$ of node $i$ in $x$.

Hypergraph beta models are families of probability distributions over
$\mathcal{H}_n$ which postulate that the hyperedges occur independently. In
details,  let $p = \{ p_e : e \in \mathcal{E}_n\}$ be a vector of probabilities whose
$e$th coordinate indicates the probability of observing the hyperedge $e$. We
will assume  $p_e \in (0,1)$. Every such vector $p$ parametrizes a
beta-hypergraph model as follows: the probability of observing the hypergraph
$x = \{ x_e, e \in E\}$ is
\begin{equation}
    \label{eq:1}
\mathbb P(x) =  \prod_{e\in E} p_e^{x_e}(1-p_e)^{1-x_e}.
\end{equation}
The graph beta model is a simple instance of this model, with  $E = \{ (i,j), 1\leq  i < j \leq n\}$.
 The ${n \choose 2}$ edge probabilities are parametrized  as
$    p_{i,j} = { e^{\beta_i + \beta_j}}/{(1 + e^{\beta_i + \beta_j})}$, for $i<j$ and  some real vector  $\beta = (\beta_1,\ldots,\beta_n)$.


Various social network modeling considerations for node interactions require a flexible class of models adaptable to those settings. Thus, we introduce three variants of the beta model for hypergraphs 
 with independent edges in the form of linear exponential families:
 \emph{beta models} for \emph{uniform hypergraphs}, for \emph{general hypergraphs}, and for \emph{layered uniform hypergraphs}. For each, 
  we provide an
exponential family parametrization in minimal form and describe the corresponding minimal sufficient statistics.
\paragraph{Uniform hypergraphs.}
The probability of a size-$k$ hyperedge $e=i_1\ldots{i_k}$ appearing in the hypergraph is parametrized by a vector $\beta\in\mathbf{R}^n$ as follows:
\begin{equation}\label{eq:un0}
p_{i_1,\ldots,i_k} = \frac{e^{\beta_{i_1}+\ldots+\beta_{i_k}}}{1+e^{\beta_{i_1}+\ldots+\beta_{i_k}}}
\end{equation}
with
$q_{i_1,
 \ldots,i_k}=1-p_{i_1,
  \ldots,i_k} = \frac{1}{1+e^{\beta_{i_1}
   +\ldots+\beta_{i_k}}}, $
for all 
 $i_1<  \cdots< i_n$.
In terms of odds ratios,
\begin{equation}\label{eq:ptobeta}
\log \frac{p_{i_1,
 \ldots,i_k}}{q_{i_1, 
  \ldots,i_k}}=\beta_{i_1}
   +\ldots+\beta_{i_k}.
\end{equation}
In order to write the model in exponential family form, we abuse notation and define for each hyperedge $e\in F$, \mbox{$\tilde{\beta}_e=\sum_{i\in e} \beta_i$.} In addition, let ${[n] \choose k}$ be the set of all subsets of size $k$ of the set $\{1,\dots,n\}$. By using (\ref{eq:1}), we obtain
\begin{equation*}
\mathbb P_{\beta}(x)  = \frac{ \exp\left\{\sum_{e\in {[n] \choose k}}\tilde{\beta}_e x_e\right\}}{\prod_{e\in {[n] \choose k}} 1+e^{\tilde{\beta}_e}}
=  \exp\left\{\sum_{i\in V}d_i(x)\beta_i-\psi(\beta)\right\},
\end{equation*}
where $d_i$ is the degree of the node $i$ in $x$. Then it is clear that the \emph{sufficient statistics for the $k-$uniform beta model} are the entries of the degree sequence vector of the hypergraph, $(d_1(x),
 \ldots,d_n(x))$, and the normalizing constant is
\begin{equation}\label{eq:unnc}
\mbox{$\psi(\beta) = \sum_{e\in {[n] \choose k}} \log(1+e^{\tilde{\beta}_e})$.}
\end{equation}

\paragraph{Layered uniform hypergraphs.}
Allowing for various size edges 
 has the advantage of controlling the propensity of each individual
to belong to a size-$k$ group independently for distinct $k$'s.
Let $r$ be the (natural bound for the) maximum
size of a hyperedge that appears in $x$. This model is then parametrized by
$r-1$ vectors in $\mathbf{R}^n$ as follows:
\[ p_{i_1,i_2,\ldots,i_k} = \frac{e^{\beta^{(k)}_{i_1}+\beta^{(k)}_{i_2}+\ldots+\beta^{(k)}_{i_k}}}{1+e^{\beta^{(k)}_{i_1}+\beta^{(k)}_{i_2}+\ldots+\beta^{(k)}_{i_k}}}
\]
where, for each $k=2,\ldots,r$, $\beta^{(k)} = (\beta_1^{(k)},\ldots,
\beta_n^{(k)})$.
There are $(r-1)n$ parameters in this parametrization. By using (\ref{eq:1}) again, we obtain
\begin{equation*}
\mathbb P_{\beta}(x)  = \prod_{k=2}^r \prod_{e\in {[n] \choose k}}\frac{e^{\tilde{\beta}^{(k)}_e x_e}}{1+e^{\tilde{\beta}^{(k)}_e}}
=  \exp\left\{\sum_{k=2}^r\sum_{i\in V}d^{(k)}_i(x)\beta^{(k)}_i-\psi(\beta)\right\},
\end{equation*}
where $d_i^{(k)}$ is the number of hyperedges of size $k$ to which node $i$ belongs in $x$. Notice that the vector of sufficient statistics in this case is
$\mathbf d =(d_1^{(2)}(x), 
\ldots, d_n^{(2)}(x), d_1^{(3)}(x), 
\ldots, d_n^{(3)}(x),\\ \ldots, d_1^{(r)}(x), 
\ldots, d_n^{(r)(x)})$,
and the normalizing constant is
\begin{equation}\label{eq:lunnc}
\mbox{$\psi(\beta) = \sum_{k=2}^r\sum_{e\in {[n] \choose k}} \log(1+e^{\tilde{\beta}^{(k)}_e})$.}
\end{equation}
\paragraph{General hypergraphs.}
In the third variant of the model we define one parameter for each node, controlling the propensity of that node to be in a relation of any size.
The probability of observing a hypergraph $x$ is thus
\begin{equation*}
\mathbb P_{\beta}(x)  =  \frac{\exp\left\{\sum_{k=2}^r\sum_{e\in {[n] \choose k}}\tilde{\beta}_e x_e\right\}}{\prod_{k=2}^r\prod_{e\in {[n] \choose k}} 1+e^{\tilde{\beta}_e}}
=  
\exp\left\{\sum_{i\in V}d_i(x)\beta_i-\psi(\beta)\right\}. 
\end{equation*}
The vector of sufficient statistics is then $\mathbf d = (d_1(x),\ldots d_n(x))$, where \mbox{$d_i(x)=\sum_{k=2}^r d^{(k)}_i(x)$}, and the normalizing constant is \mbox{$\psi(\beta) = \sum_{k=2}^r\sum_{e\in {[n] \choose k}} \log(1+e^{\tilde{\beta}_e})$.}

\section{Parameter estimation}\label{sec: parameter-estimation}
\paragraph{Iterative proportional scaling algorithms.}

 From the theory of exponential families, it is known that the MLE $\hat{\beta}$ satisfies the following system of equations:
  \begin{equation}\label{eq:mlegen}
\frac{\partial \psi(\hat{\beta})}{\partial \hat{\beta}_i}=\bar{d}_i,\nn\nn\text{for}\nn i\in\{1,\dots,n\},
\end{equation}
where $\bar{d}$ is the average observed degree sequence. By using (\ref{eq:unnc}), we then obtain
  \begin{equation}\label{eq:mlebet}
\sum_{s\in {[n]\setminus\{i\} \choose k-1}}\frac{e^{\hat{\tilde{\beta}}_s+\hat{\beta}_i}}{1+e^{\hat{\tilde{\beta}}_s+\hat{\beta}_i}}=\bar{d}_i,\nn\nn\text{for}\nn i\in\{1,\dots,n\},
\end{equation}
which is itself equivalent to
$\sum_{s\in {[n]\setminus\{i\} \choose k-1}} \hat{p}_{s,i}=\bar{d}_i$, for $i\in\{1,\dots,n\}$.

Iterative proportional scaling (IPS) algorithms fit the necessary margins of a provided table, whose elements correspond to the mean-value parameters (in this case probabilities of observing an edge). We design the following IPS algorithm for computing $\hat{p}$.

\begin{alg}\label{alg:IPS}
Define $A=(a_{i_1,\dots,i_k})$ to be an $n\times \dots\times n$  $k$-way table with margins $\bar{d}_1,\dots,\bar{d}_n$ for all its layers. Set the following structural zeros on the table: $a_{i_1,\dots,i_k}=0$ if $i_a=i_b$ for at least one pair $a\neq b$, $1\leq a,b\leq k$. (Note that there are $n(n-1)\dots(n-(k-1))$ non-zero elements in the table.)
Place $2\bar{e}/(n(n-1)\dots(n-(k-1)))$ on all other elements of the matrix, where $2\bar{e}=\sum_{i=1}^n\bar{d_i}$. Then apply the following iterative $(t+1)$st step  for every element $a_{i_1,\dots,i_k}$:
$a_{i_1,\dots,i_k}^{(t+1)}=a_{i_1,\dots,i_k}^{(t)}(F_{i_1}^{(t)}\dots F_{i_k}^{(t)})^{1/k},$
where $F_{i_b}(t)=d_{i_b}/\sum_{s\in {[n]\setminus\{i_b\} \choose k-1}}a^{(s)}_{i_b,i_s}$.
\end{alg}
IPS algorithms are known to converge to elements of the limiting matrix ($\hat{p}_{i_1,\dots,i_k}$) which are unique and preserve all the marginals (see e.g.\ \cite{bis75}). Solving the system (\ref{eq:ptobeta}) for every $1\leq i_1<\dots<i_k\leq n$ provides $\hat\beta$. Algorithm~\ref{alg:IPS} 
can be adjusted for layered uniform and general hypergraph beta models. 

For  layered $k$-uniform hypergraphs, by using (\ref{eq:mlegen}) and (\ref{eq:lunnc}) we obtain for $i\in\{1,\dots,n\}$ and $k\in\{2,\dots,r\}$,
\begin{equation}\label{eq:mlelay}
\sum_{s\in {[n]\setminus\{i\} \choose k-1}} \hat{p}_{s,i}=\bar{d}^{(k)}_i.
\end{equation}
Therefore, 
 we can apply Algorithm~\ref{alg:IPS} 
 to $(r-1)$ $k$-way tables similar to those of the $k$-uniform case, where $k$ ranges from $2$ to $r$. 
 
For general hypergraphs, we similarly obtain
\begin{equation}\label{eq:mlegenh}
\sum_{k=2}^r\sum_{s\in {[n]\setminus\{i\} \choose k-1}} \hat{p}_{s,i}=\bar{d}_i,\nn\nn\text{for}\nn i\in\{1,\dots,n\}.
\end{equation}
 In this case we apply the IPS algorithm to the following table: Define $A=(a_{i_1,\dots,i_k})$ to be a $k$-way table of size $(n+1)\times (n+1)\times\dots\times (n+1)$ consisting of labels $(\varnothing,1,2,\dots,n)$ with margins $\bar{d}_{\varnothing},\bar{d}_1,\dots,\bar{d}_n$ for all its layers, where $\bar{d}_{\varnothing}$ does not need to be known or calculated. We also set the following structural zeros in the table: $a_{i_1,\dots,i_k}=0$ if (1) $i_a=i_b\neq \varnothing$ for at least one pair $a\neq b$, $1\leq a,b\leq k$; (2) $i_1=\dots=i_k=\varnothing$ except possibly for one $i_b$.
We apply Algorithm~\ref{alg:IPS} 
 as in the $k$-uniform case except the fact that we do not fit the $\bar{d}_{\varnothing}$ margins.  We read the elements of the limiting matrix of from, $\hat{p}_{\varnothing,s}$ as $\hat{p}_{s}$, which corresponds to a lower dimensional probability.
\paragraph{Fixed Point Algorithms.}
An alternative method for computing MLE is based on 
 \cite{ChatterjeeDS2011BetaModel}.
In the $k$-uniform case, for \mbox{$i\in\{1,\dots,n\}$}, Equation (\ref{eq:mlebet}) can be rewritten as
\begin{equation}\label{eq:mlefp}
\hat{\beta}_i=\log d_i-\log \sum_{s\in {[n]\setminus\{i\} \choose k-1}}\frac{e^{\hat{\tilde{\beta}}_s}}{1+e^{\hat{\tilde{\beta}}_s+\hat{\beta}_i}}:=\varphi_i\left(\hat{\beta}\right).
\end{equation}
Therefore, in order to find $\hat{\beta}$, it is sufficient to find the fixed point of the function $\varphi$. 
\begin{alg}\label{alg:2}
Start from any $\hat{\beta}_{(0)}$ and define $\hat{\beta}_{(l+1)}=\varphi(\hat{\beta}_{(l)})$ for $l=0,1,2,\dots$.
\end{alg}
\begin{thm}\label{thm:fixed-point-alg}
If the MLE exists, Algorithm \ref{alg:2} converges geometrically fast; if the MLE does not exist there is a diverging subsequence in $\{\hat{\beta}_{(i)}\}$.
\end{thm}

The proof is omitted due to space limitations.
For the other models, the above theory can be easily generalized. For the layered models and general hypergraph models, we apply the same algorithm to obtain the fixed points of  the following functions respectively for $i\in\{1,\dots,n\}$ and $k\in\{2,\dots,r\}$ and $i\in\{1,\dots,n\}$.
\begin{eqnarray}
\varphi_i(\hat{\beta}^{(k)})&:=&\log d^{(k)}_i-\log\sum_{s\in {[n]\setminus\{i\} \choose k-1}}\frac{e^{\hat{\beta}^{(k)}_s}}{1+e^{\hat{\beta}^{(k)}_s+\hat{\beta}^{(k)}_i}}; \\ \label{eq:mlefplay}
\varphi_i(\hat{\beta})&:=&\log d_i-\log \sum_{k=2}^r\sum_{s\in {[n]\setminus\{i\} \choose k-1}}\frac{e^{\hat{\tilde{\beta}}_s}}{1+e^{\hat{\tilde{\beta}}_s+\hat{\beta}_i}}. \label{eq:mlegen}
\end{eqnarray}
\section{Simulations and Analysis}\label{sec:simulations}
\paragraph{MLE.} We use the fixed point algorithm to estimate the natural parameters for hypergraph beta models, examine non-existence of MLE and compare the layered and general variants of the model on simulated data.
Note that most dense hypergraphs, when reduced to binary relations give the complete graph, for which the MLE does not exist. In contrast, MLE is expected to exist for the hypergraph beta model in this case.  
\begin{example}\label{eg:3-uniform-beta}
We simulate a hypergraph $H=(V,F)$ drawn from the beta model for $3$-uniform hypergraphs on $10$ vertices with
$\mathbf \beta =
(-5.05,-0.57, 2.87, 4.85, 1.98, -6.69, -3.95, 5.97,\\-6.61,-4.24)$. The average simulated degree sequence of hypergraphs drawn from this model is
${\bar{d}}=(6.28,10.70,17.59,20.81,16.55,4.41,7.47,23.02,4.50, 7.17),$ and the average simulated density of the corresponding hypergraph is 0.33. Algorithm \ref{alg:2} provides the following MLE estimate using $\bar d$ as the sufficient statistic:
$\hat\beta=( -4.94,-0.58,2.81,4.76,1.94,-6.55,-3.86,5.86,\\-6.48,-4.15)$. Note that $||\mathbf \beta-\hat\beta||_{\infty}=0.14$.

For a larger example,  we select  a $\beta$  value giving rise to $3$-uniform hypergraphs on $100$ vertices with density $0.44$, and obtain a closer estimate: $||\mathbf \beta-\hat\beta||_{\infty}=0.12$.
\end{example}
\begin{example}
Theorem~\ref{thm:fixed-point-alg} guarantees that if $\hat{\beta}$ is the solution to the ML equations (\ref{eq:mlefp}), (\ref{eq:mlefplay}), or (\ref{eq:mlegen}), then the sequence of $\beta$-estimates that the fixed point algorithm produces will converge to $\hat{\beta}$; else there will be a divergent subsequence.
To detect a divergent sequence in practice, we either look for a periodic subsequence, or for a number with large absolute value in the sequence that seems to be growing, sometimes quite slowly. From \eqref{eq:un0}, since $e^{\beta_{i_k}}/(1+e^{\beta_{i_k}})$ converges to $1$ quickly ($e^{10}/(1+e^{10})>0.9999$),
for graphs with small number of nodes (i.e. far from the asymptotic behavior),
it is plausible to conclude that the corresponding mean value parameter is approximately $0$ or $1$, and hence the MLE does not exist.
Figure~\ref{fig:chainex} demonstrates MLE existence against edge densities for random hypergraphs with a fixed edge-density. Interestingly, in this restricted class, our simulations give evidence of a transition from non-existence of the MLE to existence as the density of the hypergraphs increases. The transition point seems to depend on both the number of vertices and the edge sizes allowed in the model. 
\begin{figure}
\vspace{-2mm}
\centering
\begin{tabular}{ccc}
\scalebox{0.3}{\includegraphics{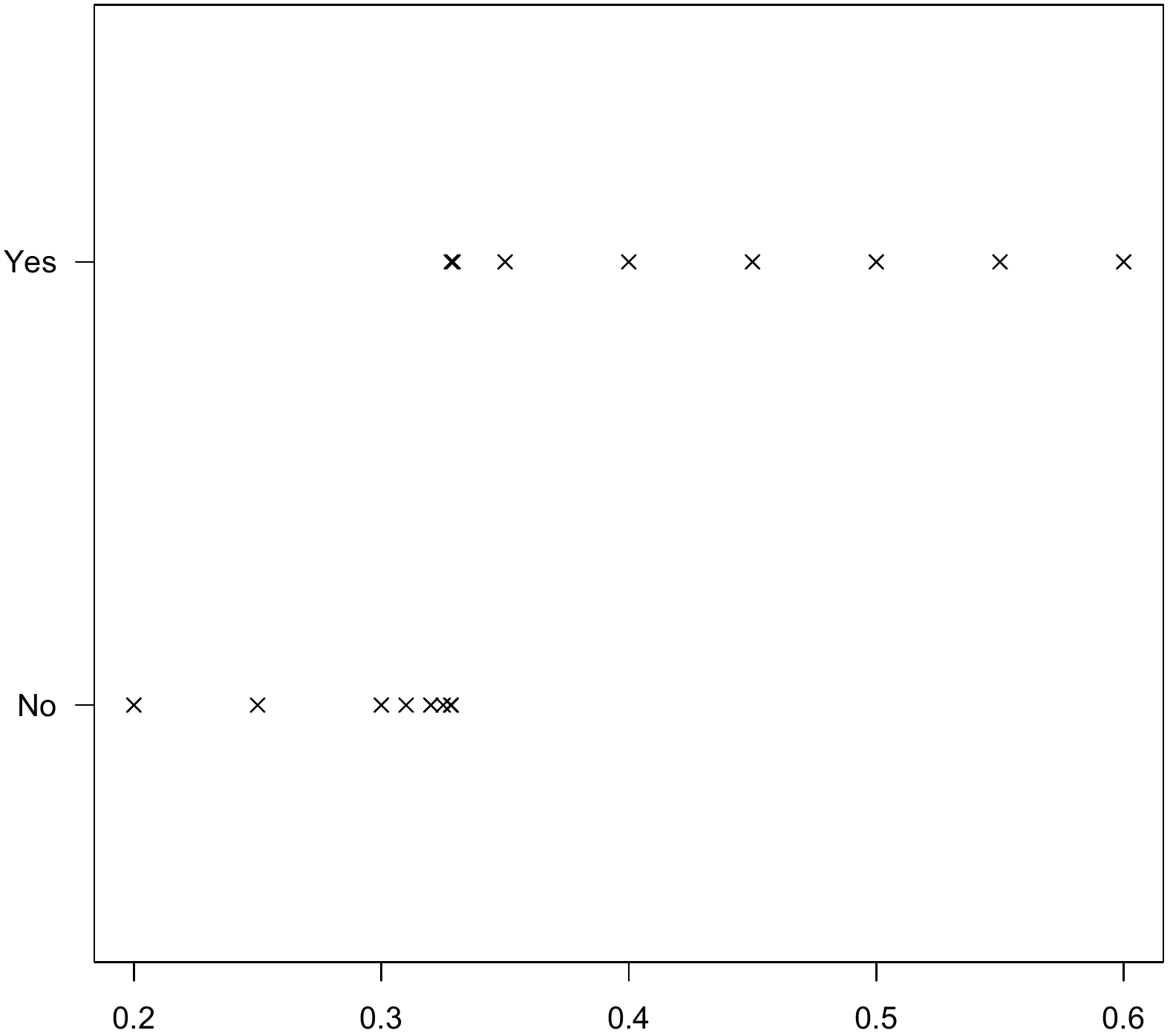}}& \scalebox{0.3}{\includegraphics{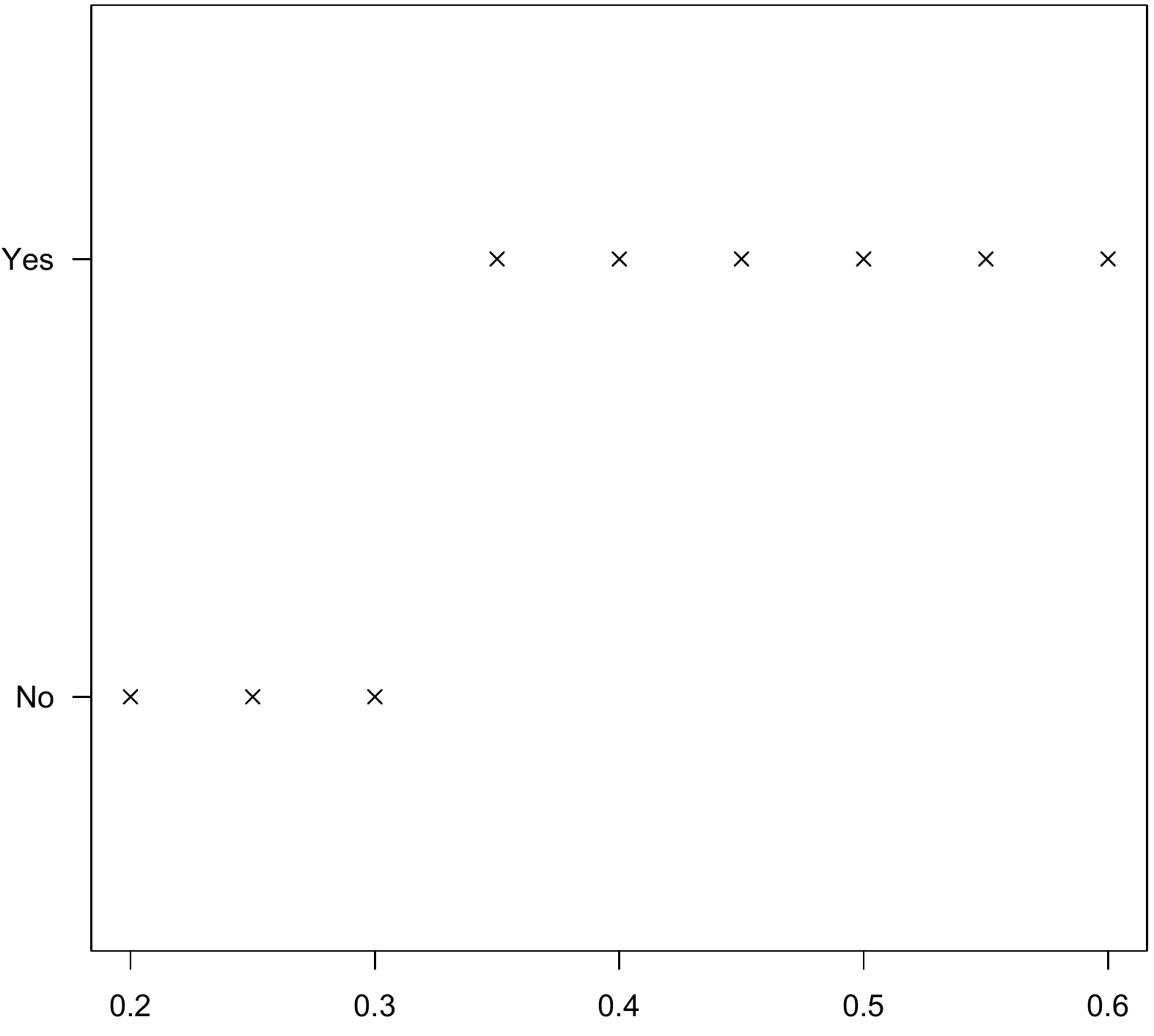}} &
\scalebox{0.3}{\includegraphics{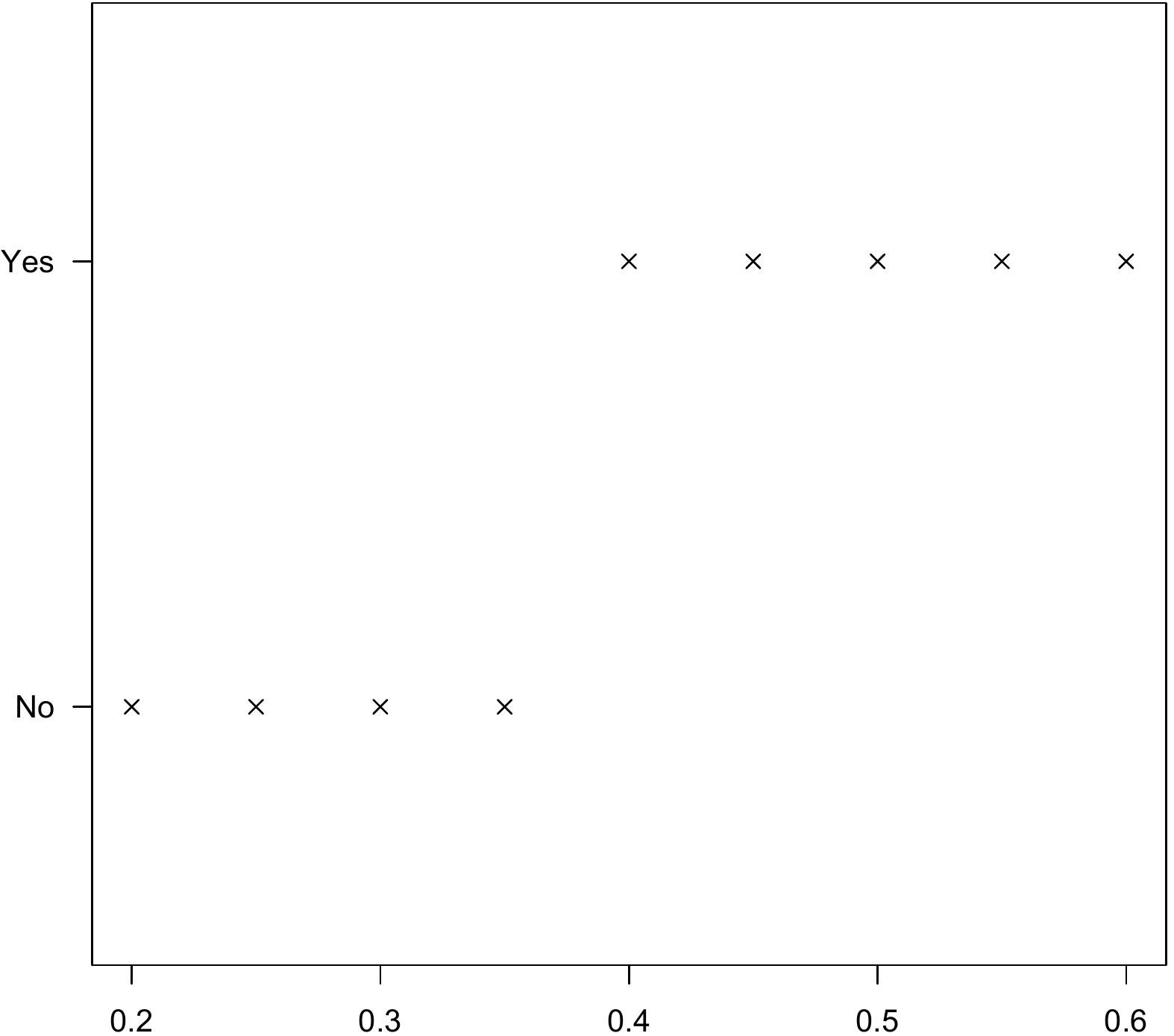}} \\ 
\scalebox{0.3}{\includegraphics{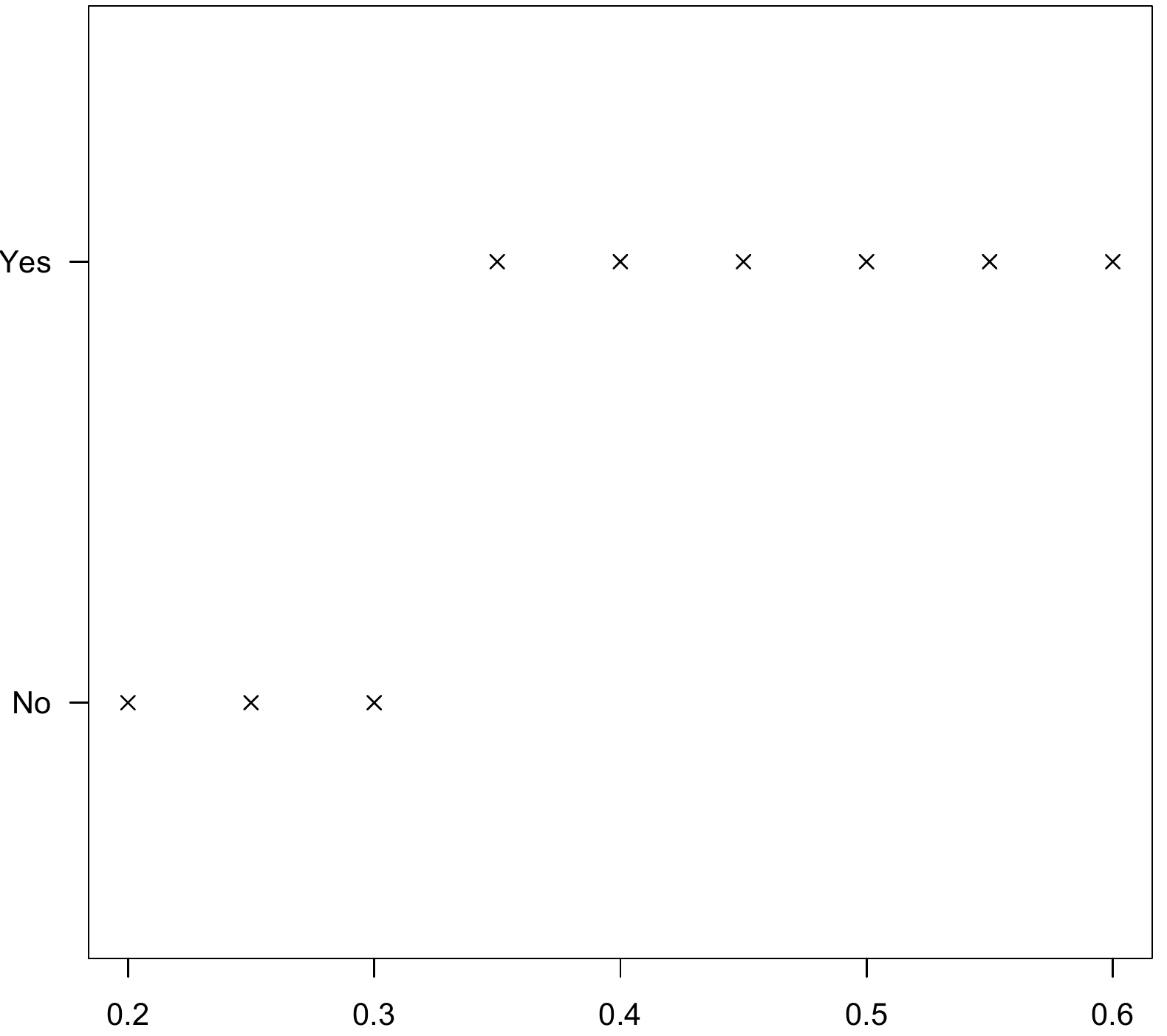}} &
 \scalebox{0.3}{\includegraphics{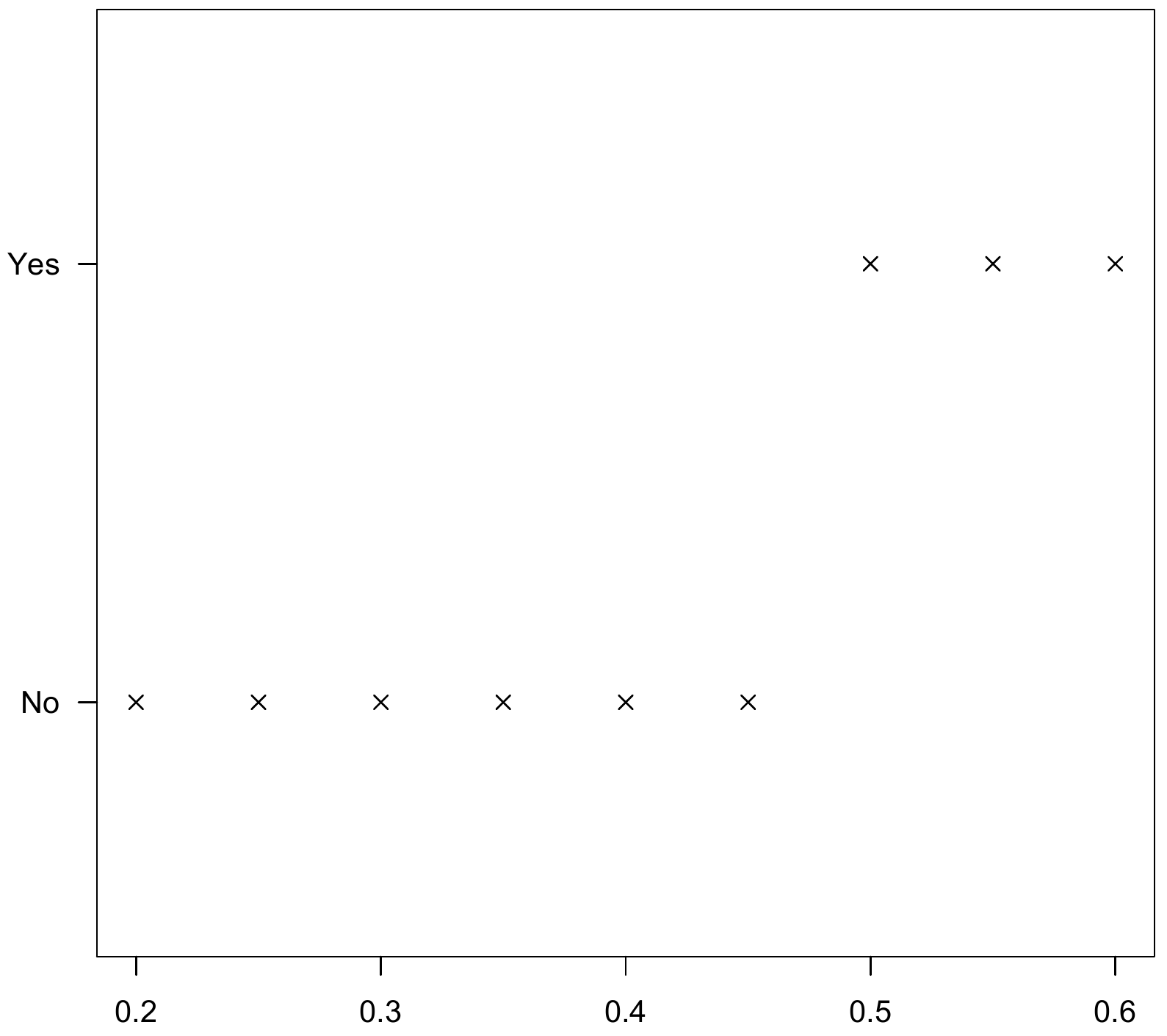}} & \scalebox{0.3}{\includegraphics{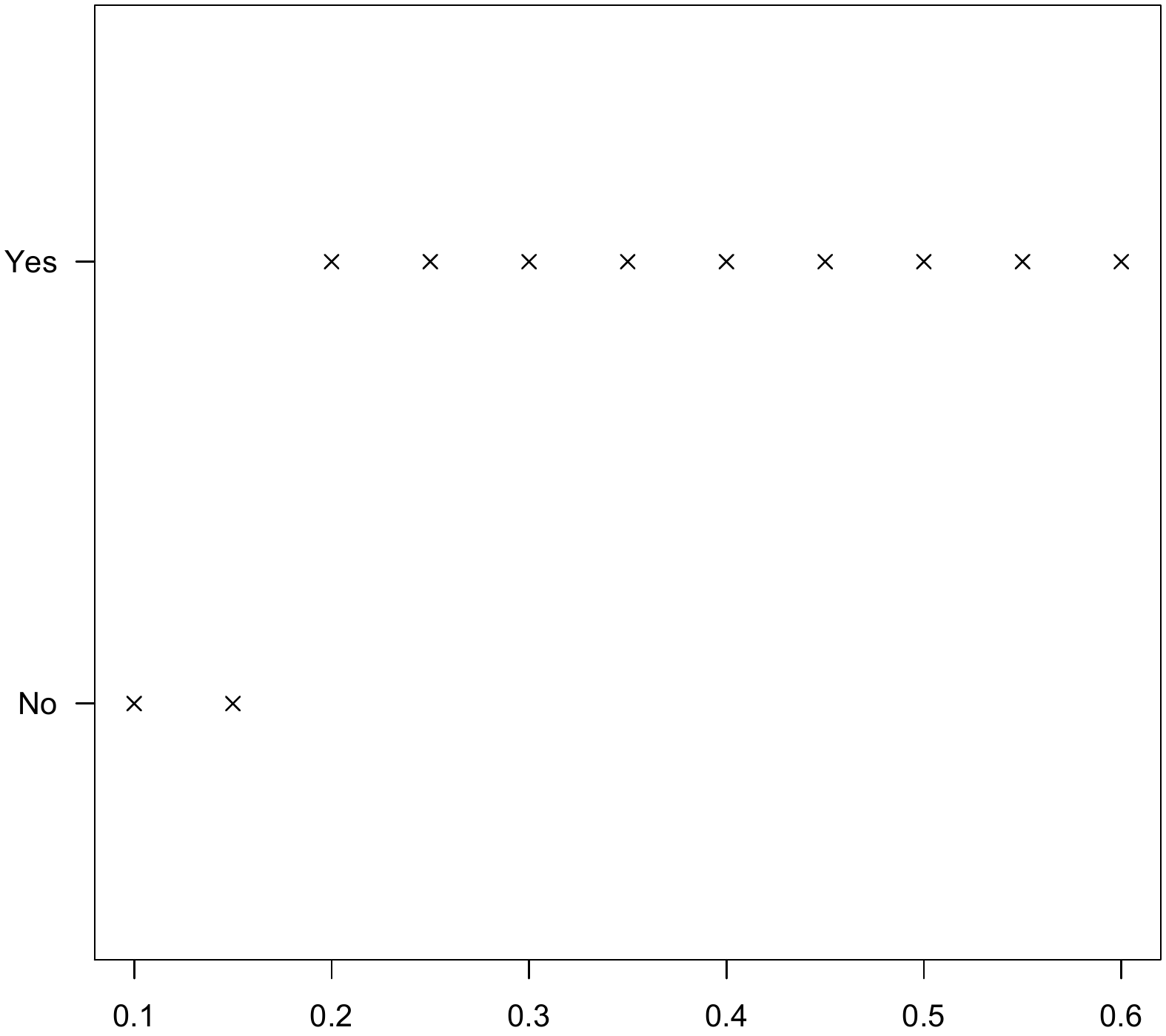}}   \\
\end{tabular}
  \caption[]{\small{MLE existence (vertical axes) for the $k$-uniform beta model on $n$ vertices against edge density (horizontal axes). Top row: $k=3$; $n=25$ (left), $n=50$ (middle), $n=100$ (right). Bottom row: $k=2$ and $n=25$ (left); $k=4$ and $n=25$ (middle); $k= \{2,3\}$ and $n=50$ (right).
  }
  }
     \vspace{-4mm}
          \label{fig:chainex}
\end{figure}
\paragraph{Model fitting: Layered versus general hypergraph beta models.}
Consider the two variants of the beta model for non-uniform hypergraphs:  the general model, with one parameter $\beta_i$ per node $i$, and the layered model, with one parameter $\beta_i^{(k)}$ per node $i$ and edge size $k$. Since the former can be considered a submodel of the latter by setting certain constraints on  $\beta_i^{(k)}$,  $k\in\{1,\dots,r\}$, we compare the fit of these two models using the likelihood ratio test with test statistics $\lambda = 2\log{\mathcal L(\hat\beta_{\text{layered}})}-2\log{\mathcal L(\hat\beta_{\text{general}})}$. Our experiments indicate that the layered model fits significantly better than the general case. Using 100 random sequences on 10 vertices, with allowed edge-sizes 2 and 3, we obtain the average observed test statistics $53.649$, in the critical region for $0.005$ significance level, $(25.188,\infty)$, for chi-square with $10$ degrees of freedom. The layered model fits significantly better for significance level 0.05 in all 100 cases, and 97 and 94 times better for significance levels 0.01 and 0.005, respectively.
\end{example}

\end{document}